\documentclass[12pt]{article}
\usepackage{amssymb}
\usepackage{amsmath,amsthm}

\newtheorem{theorem}{Theorem}

\newtheorem{proposition}[theorem]{Proposition}

\begin{document}

\title{\textbf{\ Ultrametric Cantor sets and Growth of Measure}}
\author{Dhurjati Prasad Datta\thanks{Corresponding author, email: dp\_datta@yahoo.com}, Santanu Raut\thanks{Chhatgurihati Seva Bhavan Sikshayatan High School, New Town, Coochbehar-736101, email: raut\_santau@yahoo.com} and Anuja Raychoudhuri \thanks{
Ananda Chandra College, Jalpaiguri-735101, India, email: anujaraychaudhuri@ymail.com}\and \ Department of Mathematics \and University of North Bengal,
Siliguri, West Bengal, India, Pin 734013}
\date{}
\maketitle

\baselineskip=17.5pt

\begin{abstract}
A class of ultrametric Cantor sets $(C, d_{u})$ introduced  recently (Raut, S and Datta, D P (2009), Fractals, 17, 45-52) is shown to enjoy some novel properties. The ultrametric $d_{u}$ is defined using the concept of {\em relative infinitesimals} and an {\em inversion } rule. The associated (infinitesimal) valuation  which turns out to be both scale and reparametrisation invariant, is identified with the Cantor function associated with a Cantor set $\tilde C$ where the relative infinitesimals are supposed to live in. These ultrametrics are both metrically as well as topologically inequivalent compared to the topology induced by the usual metric. Every point of the original Cantor set $C$ is identified with the closure of the set of gaps  of $\tilde C$. The increments on such an ultrametric space is accomplished by following the inversion rule. As a consequence, Cantor functions are reinterpreted as locally constant functions on these extended ultrametric spaces. An interesting phenomenon, called {\em growth of measure}, is studied on such an ultrametric space. Using the reparametrisation invariance of the valuation it is shown how the scale factors of a Lebesgue measure zero Cantor set might get {\em deformed} leading to a {\em deformed} Cantor set with a positive measure. The definition of a new {\em valuated exponent} is introduced which is shown to yield the fatness exponent in the case of a positive measure (fat) Cantor set. However, the valuated exponent can also be used to distinguish  Cantor sets with identical Hausdorff dimension and thickness. A class of Cantor sets with Hausdorff dimension $\log_3 2$ and thickness 1 are constructed explicitly.
\end{abstract}

\begin{center} 
Measure
\end{center}

\begin{center} 
\end{center}

\begin{center}
 {\em p-adic Numbers, Ultrametric Analysis and Applications}, {\bf 3}, 7-22, (2011).
\end{center}
\newpage 

\section{Introduction}

A Cantor set is a compact, perfect, totally disconnected, metrisable topological space. In this work we consider a Cantor set that is realized as a (proper) subset of  the real line. It is of measure zero if the Lebesgue measure of the
set is zero. Otherwise this has a positive measure. Such a set is also said
to be a fat Cantor set. Hausdorff measure and Hausdorff dimension are
generally considered to reveal the characteristic features of a Cantor set.
A set $C$ is said to be an $s$- set if the corresponding Hausdorff measure has a finite
non-zero value viz; $0<$ $H^{s}(C)<\infty $  \cite{falc}. Cantor sets are known to carry
a natural ultrametric structure. Here we study a family of Cantor sets with
same Hausdorff dimension. We endow a  Cantor set with a scale
invariant ultrametric structure which coincides with the usual ultrametric
only for a special choice of the valuation. For a general valuation this
offers one with an insight how a measure zero set might grow into a fat
Cantor set because of a distortion (deformation) of scales in an
infinitesimal element of the original thin Cantor set. Recently, there has 
been some interest in ultrametric Cantor sets \cite{pb, b1, sd1, sd2, ad}. Let us recall that a metric on a Cantor set is {\em regular} if the topology induced by the metric is equivalent to the usual topology. The relationship of the present class of  ultrametrics with those considered in  the so-called Michon's correspondence with rooted trees \cite{pb} will be considered separately. In the present work, we show that the valued measure induced by  our non - archimedean valuation retrieves the nontrivial measure of  \cite{pb} in a natural way. This valuation also gives a new handle to distinguish two sets with identical
Hausdorff dimension. Various other measures are available in literature to
characterise different aspects of a Cantor set. Thickness is used to study arithmetic sum and difference of Cantor sets that arise
in number theory and dynamical systems \cite{hall, newh}. However, thickness of a positive measure set is
infinite and fails to reveal the intrinsic geometric structures of such a
set. A commonly used measure is the fatness exponent \cite{far1, far2}. An analogous but, nevertheless, a  slightly
different parameter is called the exterior dimension or the uncertainty
exponent \cite{ott}. We compare our renormalised valuated parameter with these parameters.

\subsection{Main Results}

In Sec. 2, we collect, for completeness,  the basic definitions and properties that are necessary for this paper. 

\vspace{.25cm}

In Sec.3, we present, in some detail, the definition(s) of a class of inequivalent ultrametrics \cite{sd1,sd2,ad} and discuss the salient properties, mostly, in the context of the classical middle third Cantor set. The definition of ultrametrics depends on the concept of relative infinitesimals which are {\em new} elements in the neighbourhood of 0 and satisfy an {\em inversion rule}. That is to say, the increments among infinitesimals as well as between an infinitesimal and a (real) point of the Cantor set are accomplished by {\em inversions}, rather than by {\em translations} that is generally considered in standard real analysis. Some of the results presented here are new (Theorem 1 and Proposition 2): (i) The relation of  the nontrivial valuation with the Cantor function, considered in Refs. \cite{sd1, sd2}, is made precise by proving that the valuation is indeed given by  an appropriate Cantor function. (ii) The multiplicative representation that exists because of the nontrivial infinitesimals and the scale invariant ultrametric for every element of the Cantor set is verified explicitly, leading to a proof that the non-archimedean absolute value $||x||=3^{-ns}, \ x\in C, \ s=\log_3 2$ precisely corresponds to the ultrametric of \cite{pb} in the context of noncommutative geometry and (iii) the analysis of convergence of sequences of the form $\epsilon^{n-nl}, \ 0<\epsilon, l<1 $. The usual limit 0 is replaced by the constant $l$ in the present ultrametrics. This establishes the metric as well as the topological inequivalence of these scale invariant ultrametrics. The ultrametrics considered in \cite{pb} are only metrically distinct.

\vspace{.25cm}

In Sec. 4, we discuss the relationship of the scale invariant differential equation (DE)

\begin{equation}\label{eq1}
x\frac{dX}{dx}=X
\end{equation}

\noindent with a Cantor set $C$. First, we reinterpret the so-called non-smooth solutions of Ref.\cite{dp1, dp2} in the context of Cantor sets, when Cantor set elements are replaced by infinitesimal copies of an {\em inverted ultrametric Cantor set}. Infinitesimals are supposed to live in such a set which is defined as the closure of the set of gaps of a Cantor set. Eq(1) is thus well defined in the neighbourhood of each $x\in C$. The non-smooth solutions of \cite{dp2} are realised as locally constant functions (LCF) $\phi(x)$ on $C$. The nontrivial ultrametric valuation that is shown to correspond to the locally constant (LC) Cantor function associated with $C$ is, therefore, raised to the class of smooth (LC) solutions of the scale free equation (1) in the double logarithmic ($\log\log x^{-1}$) scale (c.f. eq.(26)) (Theorem 2). 

\vspace{.25cm}

In Sec.5, we study the interesting phenomenon of {\em growth of measure} that becomes meaningful in the context of the present ultrametric. (i) Starting from the observation of the {\em reparametrisation} invariance of the definition of a LC valuation in an ultrametric space, we demonstrate explicitly how a possible infinitesimal scale variation  in the reparametrisation invariant infinitesimal can lead to a {\em deformation} of scale factors of a (Lebesgue) measure zero Cantor set, so that the associated {\em deformed} Cantor set may have a positive measure (Theorem 3).  Next, (ii) we define a higher (double logarithmic) order valuation and an {\em valuated exponent} that is shown to equal the fatness exponent that is used in literature \cite{far1, far2} to characterise Cantor sets with positive measure. However, (iii) the valuated exponent is identified here with the inverse of the Hausdorff dimension of a residual  (measure zero) Cantor set that would arise as local fine structures in the neighbourhoods of the points of the original (positive measure) Cantor set. (iv) The valuated exponent is also shown to characterise a family of measure zero sets all having identical Hausdorff dimension and thickness. This is demonstrated by constructing a family of such sets with Hausdorff dimension $\log_3 2$.

\vspace{.25cm}

In the  concluding Sec.6, we remark on the implications of the inequivalent ultrametrics and the associated multiplicative structure, called {\em the generalized Euler's fatorisation} on the real number system. The nature of the asymptotic limit $x\rightarrow 0^+$ is altered significantly in the presence of relative infinitesimals.
  
\section{Basic Definitions}

A Cantor set $C$\ is  defined as a countable intersection of
finite unions of closed (and bounded) subsets of $R$. For definiteness, let $%
C\subset I=$\ $[0,1]$. \ Then, by definition, $C=$ $\overset{\infty }{%
\underset{1}{\cap }}$ $\ F_{n}=$ $\underset{n=1}{\overset{\infty }{\cap }}%
\underset{m=1}{\overset{p ^{n}}{\cup }}F_{nm}$\ where $F_{nm}\subset I$ are
closed with $F_{00}=I.$ Equivalently, $C$\ is also defined as $C=I-$ $%
\underset{i=1}{\overset{\infty }{\cup }}O_{i}$ where $O_{i}$ are open
intervals which are deleted recursively from $I$. Consequently, a Cantor set
is often defined as the limit set of an iterated function system (IFS) $f=\{f_{i}$%
\ $\mid f_i$ $:I\rightarrow I,$ $i=1,2....,p \}$ so that $C=f$\ $(C)$. For
definiteness, we consider binary Cantor sets in which each application of the
IFS removes an open interval from a closed subinterval
splitting it into two disjoint closed subintervals of the form

\begin{equation}\label{eq2}
F=F_{0}\cup O\cup F_{1}.
\end{equation}

The deleted interval $O$\ is called the gap and the two closed components
are the bridges. As an example let us consider a middle $\alpha $ Cantor set 
$C_{\alpha }$ which arises as the limit set under the IFS 
\begin{equation}\label{eq3}
f_{i}(x)=\beta x+i(1-\beta ),\text{ \ }i=0,1
\end{equation}

\noindent where the scale factor $\beta $ is defined by $\alpha +2\beta =1$.
Each iteration of the IFS  removes an open interval (i.e. a gap)
of length proportional to $\alpha $ from a closed subinterval  of $I$, leaving out two bridges of size proportional to $\beta $ each. The IFS (\ref{eq2}) satisfies
the {\em open set condition} (OSC)  if there exists a non-empty bounded open set $S$ such
that $\underset{i}{\cup }$ $f_i$ $(S)\subseteq S.$ It follows accordingly
that $\beta \in (0,\frac{1}{2}).$ Since the total length of the disjoint
open intervals viz., $\overset{\infty }{\underset{i=1}{\Sigma }}$ \ $\mid
O_{i}\mid $ $=$\ \ $\overset{\infty }{\underset{i=1}{\Sigma }}$ $\alpha
(2\beta )^{i-1}=1,$ the middle $\alpha $ Cantor set is of measure zero with
the Hausdorff dimension $s=\frac{\log 2}{\log \frac{1}{\beta }}$. More
generally, when $q$ open intervals each of size $\alpha $ are deleted
leaving out $p$ equal closed intervals of size $\beta $ so that $q\alpha
+p\beta =1$, then the OSC gives $\beta \in (0,\frac{1}{p}).$ The length
of the deleted open intervals add up to $1$ viz., $\Sigma (q\alpha )(p\beta
)^{n-1}=1.$ The corresponding measure zero set $C_{\alpha ,p}$ has the
Hausdorff dimension $\frac{\log p}{\log \frac{1}{\beta }}.$

\vspace{.25cm}

Returning to the discussion of the binary Cantor set we recall that the set $%
C_{\alpha }$\ is also a homogeneous and uniform Cantor set. It is
homogeneous since the scale factors in each component of the IFS are same. The set is uniform because each deleted open
interval also is of constant proportion $\alpha $ of the length of the
previous (defining) closed interval.

\vspace{.25cm}

A positive 1-set $\tilde{C}$, on the other hand, is obtained if the deletion
process removes open intervals of variable sizes.

\vspace{.25cm}

\textbf{Example 1:} Let at each step we remove $\alpha _{n}$\ \ portion of
the length of each component of the previous closed set $F_{n-1}$ so that$\
\ F_{n-1}=F_{n0}\cup O_{n}\cup F_{n1}$ and $\mid O_{n}\mid =\alpha _{n}\mid
F_{n-1}\mid ,\mid F_{n0}\mid $ $=$ $\mid F_{n1}\mid =\frac{1}{2}(1-\alpha
_{n})\mid F_{n-1}\mid$. By induction, each of $2^{n}$\ \ components of $F_{n}
$\ has length $\mid F_{ni}\mid =\frac{1}{2^{n}}\underset{0}{\overset{n}{\prod }
}(1-\alpha _{j}),$ $\ i=1,2,\cdots 2^{n}.$ Consequently, $m$ $(\tilde{C})=$ $
\underset{n\rightarrow \infty }{\lim }\mid F_{n-1}\mid $ $=$ $\underset{0}{
\overset{\infty}{\prod}}(1-\alpha _{i})>0$ when $\Sigma $ $\alpha
_{n}<\infty .$

\vspace{.25cm}

\textbf{Example 2:} Suppose at the $n$th step an open interval of length $%
\frac{\delta }{3^{n}},$ $(0<\delta <1)$\ is removed from each of the $2^{n}$
components of $F_{n}$. The length of each component of $F_{n}$ is $\frac{1}{%
2^{n}}$ $(1-\frac{\delta }{3}-\cdots \frac{2^{n-1}\delta }{3^{n}}).$ The sum
of the lengths of all the open intervals removed is $\Sigma $\ $\frac{%
2^{n}\delta }{3^{n+1}}=\delta $, so that $m$ $(\tilde{C})=1-\delta.$

\vspace{.25cm}

In both the above examples, $0<$ $H^{s}(\tilde C)=m(\tilde C) = l<\infty $
when $s=1$. Even as the Cantor set 
$\tilde{C}$\ is nowhere dense in $I$, it has a non-zero measure because of
the fact that the defining iteration process now removes the open interval
at a slower rate in comparison to the middle $\alpha $ set $C_{\alpha }$. In
fact, the relative difference of the lengths of the deleted open middle $%
\alpha $ and $\alpha _{n}$ intervals $O_{\alpha n}$ and $\tilde{O}_{\alpha n}
$ respectively is given by

$$\bigg | \mid O_{\alpha n}\mid -\mid  \tilde{O}_{\alpha n} \mid
\bigg | = \bigg | \alpha \beta ^{n}- \frac{1}{2^{n}} \alpha _{n}\underset{%
1}{\overset{n}{\Pi }}(1-\alpha _{i})\bigg |  = \bigg | 1-\frac{\alpha _{n}}{%
\alpha }\underset{1}{\overset{n}{\Pi }}(\frac{1-\alpha _{i}}{1-\alpha })\bigg |
 | O_{\alpha n} |  \geq   | 1-\gamma  | | O_{\alpha n}| $$

\noindent where $\frac{\alpha _{n}}{\alpha }\underset{1}{\overset{n}{\Pi }}(\frac{%
1-\alpha _{i}}{1-\alpha })$ $\rightarrow $ $\gamma $ for $n\rightarrow
\infty.$ Note that the above lower bound exists, otherwise $\tilde{C}$ would
have been a set of measure zero. For the modified middle $\frac{1}{3}$rd set
(Example 2) one has the exact equality

$$\bigg | \mid O_{\frac{1}{3}n}\mid -\mid  \tilde{O}_{\delta } \mid \bigg | 
=(1-\delta )\mid O_{\frac{1}{3}n}\mid $$

 The emergence of the positive measure of $\tilde{C}$ can, therefore,
be explained in a dynamical sense provided the said set $\tilde C$ is seen
as being  evolved from a given zero measure set because of a possible {\em principle}
allowing for a deformation of scales in the deletion process. As is evident the
relative scaling of the lengths of infinitesimal elements of the deleted
open intervals $\tilde{O}_{\alpha n}$ over the corresponding $\alpha $
interval $O_{\alpha n}$, indeed captures the origin of the positive measure
even in the totally disconnected perfect set $\tilde{C}$. In the following
we offer a new ultrametric explanation of the growth of the positive measure
of $\tilde{C}$ over $C_{\alpha }$. The class of ultrametrics that we
consider is not only {\em scale invariant} (in the sense of a power law), but also {\em reparametrisation invariant} (that is, the invariance under a reparatrisation of the form $X(t)\rightarrow \tilde X(t)=X(f(t))$ where the otherwise arbitrary function $f$ satisfies the conditions $f^{\prime}>0$, along with the boundary condition $ f(0)=0, \ f(1)=1$).
As will be explained later, the scale variation can, therefore, be interpreted  as a reflection of the underlying reparametrisation invariance of the valuation. Incidentally, we note that for a measure 1 set the functions defining an  IFS may not have a closed form\cite{phyd}.

\subsection{Thickness}

The measure of thickness of a Cantor set has various applications in number theory and
dynamical systems. To recall the definition, let $F_{i}=I$\ $-\underset{l=0}{%
\overset{i}{\cup }}O_{l}$\ form a defining sequence of the Cantor set $C$.
The $2^{i}$ components of $F_{i}$ are the closed intervals $F_{ij},$ $%
j=1,2,\cdots 2^{i},$\  which are the bridges. The deleted intervals $O_{i}$
are the gaps of $C$. Let $O_{i_k}$ denote the open deleted  subinterval of a bridge $F_{ij}$ of $F_{i}$ which divides $F_{ij}$ into two smaller bridges $F_{ij}^{L}$ and $F_{ij}^{R}$ of $F_{i}$. Let

$$\tau  (F_{i})= \underset{j}{\inf } \left\{ \frac{\mid F_{ij}^{L}\mid 
}{\mid O_{i_k}\mid },\frac{\mid F_{ij}^{R}\mid }{\mid O_{i_k}\mid }\right\}$$

The thickness $\tau (C)$ is defined by $\tau (C)=$ $\underset{i}{\sup }$ $%
\tau $ $(F_{i})$, where $\sup $ is evaluated over the defining sequences of $%
C.$ For a set $A$ containing an interval, $\tau (A)=$\ $\infty ,$\ by
definition.

\vspace{.25cm}

For the middle $\alpha $- Cantor set $C_{\alpha }$, it follows that

$$\frac{\mid F_{ij}^{L}\mid }{\mid O_{i_k}\mid }=\frac{\mid F_{ij}^{R}\mid }{
\mid O_{i_k}\mid }=\frac{\beta \mid F_{ij}\mid }{\alpha \mid F_{ij}\mid }=
\frac{\beta }{\alpha }$$

\noindent so that $\tau (C)=$ $\frac{\beta }{\alpha }$.

\vspace{.25cm}

For a positive measure set $\tilde{C}$, on the other hand, one has $\frac{
\mid F_{ij}^{L}\mid }{\mid O_{i_k}\mid }=\frac{\mid F_{ij}^{R}\mid }{\mid
O_{i_k}\mid }=\frac{1-\alpha _{n}}{2\alpha _{n}}$, leading to $\tau (\tilde{C})
$ $=\infty ,$ as expected.

\subsection{Cantor Function} 
A Cantor function is a  nonconstant and non-decreasing continuous function $\phi :[0,1]\rightarrow [0,1]$ such that $\phi^{\prime}(x)=0$ a.e. with points of nondifferentiability  $x$ lying, for instance, in the Cantor set $ C_{\alpha}$ . 
\vspace{.25cm}
 
 To construct $\phi$ explicitly, let $\phi (0)=0,\ \phi (1)=1.$ Assign $\phi (x)$ a constant value $\phi (x)=i2^{-n},\ i=1,2,....,2^n-1$ on each of 
the deleted open intervals (including the end points of the deleted
interval) of $C_{\alpha}$. 
Next, let $x$ $\in
C_{\alpha}$. Then, at the $n$th iteration,  $x$ belongs to the interior of exactly one of the $2^{n}$
remaining closed intervals each of length $\beta^n $. Let $[a_n, b_n]$
be one such intervals. Then $b_n-a_n=\beta^n$. Moreover,  $\phi(b_n)-\phi(a_n)=2^{-n}$.
At the next iteration,
assuming $x\in [a_{n+1}, b_{n+1}]$, 
( $a_n=a_{n+1}$), say, we have $\phi (a _{n})\leq  \phi (a_{n+1})<\phi(b _{n+1})\leq \phi (b _{n})$. 
Define $\phi (x)=\underset{n\rightarrow \infty }{\lim }\phi (a _{n})=\underset{n\rightarrow \infty}{\lim }\phi (b _{n})$. Then $\phi :[0,1]\rightarrow [ 0,1]$
is a continuous, non-decreasing function. Also $\phi ^{\prime }(x)=0$ for $%
x\in $ $I\backslash C_{\alpha}$ when it is not differentiable at any \ $x\in C_{\alpha}.$ (c.f., \cite{sd1, sd2}).

\subsection{Ultrametric}

The topology of a Cantor set is equivalent to an ultrametric topology. A
point $x\in C_{\alpha }$ has the unique infinite word representation

$$x=(1-\beta )\underset{0}{\overset{\infty }{\Sigma }}x_{i}\beta ^{i} 
=x_{0}x_{1}x_{2}\cdots, \ x_{i}\in \{0,1\} $$

Let $L(x,y)=n$\, such that $x_{0}=y_{0,\cdots ,}x_{n-1}=y_{n-1},$ $x_{n}\neq
y_{n},$ $x,y\in C_{\alpha }.$
The ultrametric $\tilde d_{u}$\ is defined by $\tilde d_{u}(x,y)=p^{-L(x,y)}$\ for any $p>1.$
This ultrametric is equivalent to the usual metric

$$C_{1}\tilde d_{u}(x,y)\leq d(x,y)\leq C_{2}\tilde d_{u}(x,y)$$

\noindent  for two positive constants 
$C_{1}$ and $C_{2}$, where $d(x,y)$ is the usual metric.

\vspace{.25cm}

The Cantor set thus consists of towers of closed balls (intervals) with
countable intersection property. Further the fundamental neighbourhood
system of any point consists of clopen balls.

\vspace{.25cm}

We now define an inequivalent class of ultrametrics, called valuated
ultrametrics, on a Cantor set. The definition of the ultrametric makes use
of a concept of relative infinitesimals introduced recently \cite{sd1, sd2, ad}.

\section{Inequivalent Ultrametrics}

Let $C$ be a (measure zero) Cantor set. We consider an {\em infinitesimal gap} $
O_{\inf }$ of $C$ in $I=[0,1].$ Given an arbitrarily small $x\in C-\{0\}$ (in the sense that  $x\rightarrow 0^{+}$\ on $C-\{0\}$), $\exists $ $\epsilon \in I$ and $\epsilon<x$ and an open interval $\tilde I\subset (0,\epsilon)$ such that  $\tilde I\cap C=\Phi $, the null set. This follows from the total disconnectedness of $C$. An element $\tilde{x}$ in $\tilde I$ satisfying $0<\tilde{x}<\epsilon <x$ and the {\em inversion} rule

\begin{equation}\label{eq4}
\frac{\tilde{x}}{\epsilon }=\lambda ^{-1}(\epsilon )\frac{\epsilon }{x}
\end{equation}

\noindent for a real constant $\lambda $ $(0\ll\lambda \leq 1)$\ is called a
{\em relative infinitesimal} relative to the {\em scale} $\epsilon $. For each choice of $x$ and $\epsilon$, we have a unique $\tilde x$ for a given $\lambda \in (0,1)$. Consequently, by varying $\lambda$ in an open subinterval of (0,1), we get an open interval of relative infinitesimals in the interval $(0,\epsilon)$, all of which are related to $x$ by the inversion formula.  The infinitesimal gap $O_{\inf}\subset \tilde I$ is, by definition, the set of these relative infinitesimals satisfying the inversion rule (\ref{eq4}), as $\epsilon\rightarrow 0$, in an asymptotic sense (c.f. Remarks 1 and  2). In the limit $\epsilon \rightarrow 0$, $O_{\inf}=\Phi$, in the usual topology. However, the corresponding set of {\em scale invariant infinitesimals} $\tilde O_{\inf}=\underset{\epsilon \rightarrow 0}{\lim} \{\tilde X |\tilde X={\frac{\tilde x}{\epsilon} \approx \mu\epsilon^{\alpha} \pm o(\epsilon^{\beta})}\}$ where $\mu$ is a constant and $1>\beta> \alpha\geq 0$, may be a non-null subset of (0,1) (for instance, when $\alpha=0$, in particular) (for an explanation of the asymptotic expansion of $\tilde X$ see Remark 2.1). Notice that constants $\alpha$, $\beta$ and $\mu$ may, however, depend on $\lambda$. Notice also that the infinitesimal gap $O_{\inf}=O(x, \epsilon, \lambda)$ depends  on  $\epsilon$, but apparently also on the arbitrarily small element $x$ of the Cantor set along with the parameter $\lambda$ appearing in the inversion law. But $x$ and $\epsilon$ are very closely related,  so that $O_{\inf}$ essentially depends only on $\epsilon$ and $\lambda(\epsilon)$.

For a point $x$ from a Cantor set $C$, it is natural to assume that the scale $\epsilon$ is determined by the privileged scale of the Cantor set. Two relative infinitesimals 
$\tilde{x}$\ and $\tilde{y}$\ must necessarily satisfy the condition $0$\ $<$
$\tilde{x}$\ $\leq $ $\tilde{y}$\ $<$ $\tilde{x}$\ $+$ $\tilde{y}$\ $
<\epsilon .$ As indicated already, the inversion rule maps an open interval of (relative) infinitesimals of size determined by the parameter $\lambda $\ to an arbitrarily small element $x$ of $C$. These
relative infinitesimals are endowed with a scale free (non-archimedean)
absolute value $\mid \cdot \mid_u :$ $\ O_{\inf }\rightarrow \lbrack 0,1]$\
defined by
\begin{equation}\label{eq5}
\mid \tilde{x}\mid_u = \underset{\epsilon \rightarrow 0}{\lim }\log
_{\epsilon ^{-1}}{\frac{\epsilon }{\mid \tilde{x}\mid }}, \  \tilde{x}%
\neq 0 
\end{equation}
\noindent and $|\tilde x|_u=0$, for $\tilde{x}=0$. 

The proof follows in two steps (see, for instance, \cite{ad}). One first verifies that $\mid
\cdot \mid _{u}$\ defined by (\ref{eq5}) is a semi-norm viz: (i) \ $\mid \tilde{x}
\mid_u $ $>0, $ $\tilde{x}\neq 0$\ \ (ii) $\mid -\tilde{x}\mid_u $\ $=$ $\ \mid 
\tilde{x}\mid_u $\ \ (iii) $\mid \tilde{x}+\tilde{y}\mid _{u}\leq \max
\left\{ \mid \tilde{x}\mid _{u},\mid \tilde{y}\mid _{u}\right\} ,$\ known as
the stronger triangle\ inequality. (For the property (ii), $\mid \cdot \mid
_{u}$\ should be defined on an infinitesimal neighbourhood of 0 in $[-1,1]$).

{\bf Remark 1}: 
As remarked already, the set of infinitesimals $O_{\inf} =\Phi$ when $\epsilon\rightarrow 0$. However, the corresponding asymptotic expression for the scale free (invariant) infinitesimals  is  nontrivial, in the sense that the associated valuations (\ref{eq5}) can be shown to exist as finite real numbers. 

Choosing $\epsilon=\beta^r$, the Cantor set scale factor,  the scale free infinitesimal gaps can be identified as $\tilde O^m_{\rm inf}= (0, \beta^m)$ when $\epsilon\rightarrow 0$ is realised as $n\rightarrow \infty$, $r=n+m, \ m=1,2,\ldots$. Assign nonzero constant valuation $|\tilde x_m|_u=\alpha_m$ $\forall\ \tilde X_m= \tilde x_m/\epsilon\in \tilde O^m _{\inf}$. The set of all possible scale free infinitesimals $\cup\tilde O^m_{\rm inf}\subset (0,1)$ is now realised as nested clopen circles $S_m: \{\tilde x_m:|\tilde x_m|_u=\alpha_m\}$.  The ordinary 0 of $C$ is  replaced by this set of scale free infinitesimals $0\rightarrow \ {\bf 0} =O_{\inf}/\sim = \{0, \cup S_m\}$, $\bf 0$ being the equivalence class under the equivalence relation $\sim$, where $x \sim y$ means $|x|_u=|y|_u$ . The existence of $\tilde x$ could also be concieved dynamically as a computational model \cite{sd1,sd2,ad}, in which a number, for instance, 0 is identified as an interval $[-\epsilon, \epsilon]$ at an accuracy level determined by  $\epsilon=\beta^n$.

{\bf Remark 2}:
1. The concept of infinitesimals and the associated absolute value considered here become significant only in a limiting problem (or process), which is reflected in the explicit presence of ``$\underset{\epsilon\rightarrow 0}{\lim}$" in the relevant definitions. Recall that for a continuous real valued function $f(x)$, the statement $\underset{x\rightarrow 0}{\lim} \ f(x)=f(0)$, means that ${x\rightarrow 0}$ essentially is $x=0$. This may be considered to be a {\em passive} evaluation (interpretation) of limit. The present approach is {\em dynamic}, in the sense that it offers not only a more refined evaluation of the limit, but also provides a clue how one may induce new (nonlinear) structures (ingredients) in the limiting (asymptotic) process. The inversion rule (\ref{eq4}) is one such nonlinear structure which may act nontrivially as one investigates more carefully the {\em motion} of a real variable $x$ (and hence of the associated scale $\epsilon < x$) as it goes to 0 more and more accurately. Notice that at any ``instant",  elements defined by inequalities $0<\tilde x<\epsilon<x$ in the limiting process, are well defined; relative infinitesimals are {\em meaningful} only in that {\em dynamic} sense (classically, these are all zero, as $x$ itself {\em is } zero). Scale invariant infinitesimals $\tilde X$, however, may or may not be zero classically. $\tilde X=\mu \ (\neq 0), $ a constant,  for instance, is nonzero even when $x$ and $\epsilon$ go to zero. On the other hand, $\tilde X= \epsilon^{\alpha}, \ 0<\alpha<1$, of course, vanish classically, but as shown below, are nontrivial in the present formalism. As a consequence, relative (and scale invariant) infinitesimals may be said to {\em exist} even as real numbers in this dynamic sense. The acompaning metric $|.|_u$, however, is an ultrametric.

2. However, a genuine (nontrivial) scale free  infinitesimal $\tilde X$  can not be a constant. Let $\tilde x_0=\mu\epsilon, \ 0<\mu <1$, $\mu$ being a constant. Then $|\tilde x_0|_u=\underset{\epsilon\rightarrow 0}{\lim} \log _{\epsilon}\mu =0$, so that $\tilde x_0$ is essentially the trivial infinitesimal 0 (more precisely, such a relative infinitesimal belongs to the equivalence class of 0).

3. The scale free infinitesimals of the form $\tilde X_m\approx \epsilon^{\alpha_m}$ go to 0 at a slower rate compared to the linear motion of the scale $\epsilon$. The associated nontrivial  absolute value $|\tilde x_m|_u=\alpha_m$ essentially quantifies this decelerated motion.

\vspace{.25cm}

Next, one notes that\ $O_{\inf }$\ is the union of countable family of
disjoint clopen balls (this is, in fact, true even for the natural ultrametric). Recall that an (a) open (closed) ball is defined, as
usual, using the seminorm \ \ $\mid \cdot \mid _{u}$ \ viz : \ a set $
B_{r}(a)=\left\{ \text{ }\tilde{x}\text{ :}\mid \tilde{x}-a\mid
_{u}<r\right\} $\ \ is open while $\bar B_{r}(a)=\left\{ \text{ }\tilde{x}\text{ :
}\mid \tilde{x}-a\mid _{u}\leq r\right\} $\ is a closed ball. We now assume
that the mapping \ \ $\mid \cdot \mid _{u}:O_{\inf }\rightarrow \lbrack 0,1]$
\ \ is constant on each of the component balls $B(a_{n})\subset O_{\inf
}\subseteq \cup B(a_{n}).$ Again for a given $\epsilon $\ the closure of $
O_{\inf }$\ \ \ \ \ is compact and so is covered by a finite number of
clopen balls $B(a_{i}),$ $i=0,1,\cdots ,n.$ Consequently, $\mid \cdot \mid_u $
\ is discretely (finitely) valued and so there exists a constant $0<\sigma_{\epsilon} <1$\ such that \bigskip

$$
\mid \tilde{x}\mid_u =\sigma_{\epsilon} ^{s_i}, \ \tilde{x}\in B(a_i), \ i=0,1,\cdots ,n
$$

\noindent for an ascending sequence of valuations $%
0<s_{0}<s_{1}<\cdots <s_{n}$. For an ascending sequence $\alpha _{i}>0,$\
the above can be written alternatively as

\begin{equation}\label{eq6}
\mid \tilde{x}\mid_u =\alpha _{i\text{ }}\sigma_{\epsilon} ^{\tilde s},\tilde{x}\in
B(a_i),i=0,1,\cdots ,n																																		
\end{equation}

\noindent where $\tilde{s}=s_{n.}$  It follows from (\ref{eq5}) that $\tilde x_i\approx \epsilon\epsilon^{\alpha_i\sigma_{\epsilon}^{\tilde s }+\mu_{\epsilon}}$, where $\mu_{\epsilon}$ goes to zero faster than $\sigma_{\epsilon}$ as $\epsilon\rightarrow 0^+$.

\vspace{.25cm}

Now, as in Remark 1, the limit $\epsilon\rightarrow 0^+$ on infinitesimal scales could be evaluated via a countable number of different sequences of the form $\epsilon^{m+n}$ as $n\rightarrow \infty$. For each such sequence, there exists a special value of $\sigma=\sigma_m$, all of which may be assumed to arrange in a descending order. As a consequence, $\sigma$ assumes values from a decreasing sequence of scales (of the form $\sigma_m$), called the {\em secondary} scales. The sequence  $\epsilon^n$ defines a set of  {\em primary} scales. 

\vspace{.25cm}

Next, the (infinitesimal) valuation $|\tilde x|_u=v(\tilde x)$ (say), as a mapping from $O_{\inf}$ to $I\subset R$ is continuous. The equation (\ref{eq6}), however, defines $v(.)$ only for points in the clopen balls $B(a_i), \ i=1,2, \cdots$. The definition can be extended continuously over the entire set $O_{\inf}$ for points outside the clopen balls.  
Indeed, let for a given primary scale $\epsilon$, $\sigma_i$ be the secondary scale. Let also that $y\in O_{\inf}/\cup B(a_{i})$. Then there exist $y_i\in B(a_i), \ y_{i+1}\in B(a_{i+1}) $ such that $y_i<y<y_{i+1}$, and $v(y_{i+1})-v(y_i)= (\alpha_{i+1}-\alpha_i) \sigma_i$. Clearly, the sequence $v(y_{i+1})$ is increasing and $v(y_i)$ is decreasing. Thus,  $v(y):=\lim v(y_i)$, as $i\rightarrow \infty$. We have thus proved that the scale invariant valuation $v(\tilde x)$ is indeed given by a Cantor function. Conversely, given a Cantor function $\phi(x), \ x\in [0,1]$, one can define a set of infinitesimals by the asymptotic formula $\tilde x\approx \epsilon\epsilon^{\phi(\tilde x/\epsilon)}$ as $\epsilon\rightarrow 0$ that is assumed to live in a nontrivial neighbourhood of 0.

\vspace{.25cm}

With this class of valuations, the
seminorm now extends to a non-archimedean absolute value, satisfying also the
product rule (iv) $\mid \tilde{x}$ $\tilde{y}\mid_u $ $=$ $\mid \tilde{x}%
\mid _{u}\mid $ $\tilde{y}\mid_u.$

\vspace{.25cm}

We have thus proved the 

\begin{theorem}
The infinitesimal valuation (\ref{eq5}) is a non-archimedean absolute value, and is given by a Cantor function associated with  the Cantor set containing the  relative infinitesimals. Conversely, given a Cantor function, there exists a class of infinitesimals that live in an extended ultrametric neighbourhood of 0.
\end{theorem}

Notice that an infinitesimal gap of the form $(0,\epsilon ),$ $\epsilon
\rightarrow 0^{+}$ is a segment of a line (i.e. an open interval) in the
usual topology inherited from the real line. In the topology induced from
the ultrametric $\mid \cdot \mid _{u},$ the interval $(0,\epsilon )$ instead
is realized as a totally disconnected set and so itself is equivalent to a
Cantor set $\tilde{C}$. The structure of $\tilde{C}$\ can, in general, be
very arbitrary and even be independent of the original Cantor set $C$. In
ref.\cite {sd1, sd2}, we have, however, shown how this new ultrametric
valuation for relative infinitesimals can lead to a topology identical to
that of the original Cantor set $C$. In the following we give a
brief out line of the procedure and also comment on the origin of the
inequivalence of this ultrametric.

\subsection{Middle Third Cantor Set: An Example}

Let $C$ be the standard middle $\frac{1}{3}$rd Cantor set. As will become
clear our discussion will apply generically to any measure zero Cantor set
considered in this paper. The Cantor set $C$ offers us with a privileged
set of scales $\epsilon _{n}=3^{-n}$. For a sufficiently large $n$ viz : as $
n\rightarrow \infty,$ suppose an infinitesimal gap of the form $(0,\epsilon _{n+m})$\
is decomposed into  a finite number $M$\ of open subintervals $\tilde I_{i}$\ of
relative infinitesimals with constant valuations defined by

\begin{equation}\label{eq7}
\mid \tilde{x}_{i}\mid _{u}=i\text{ }2^{-n}\text{ , }i=1,2,\cdots ,M, \ \tilde{x_i}\in I_{i}.
\end{equation}

\noindent The valuation assigned by (\ref{eq7}) is the triadic
Cantor function $\phi :I$\ $\rightarrow I$\ so that $M=2^{m}-1$\ corresponding to the scale $\epsilon _{m}=3^{-m}$. We call infinitesimals $\tilde{x}$ leaving in the island $\tilde I_{i}\subset (0,\epsilon _{n+m})$\  \emph{valued infinitesimals} having the valuation (\ref{eq7}) induced by the Cantor function. Any element $x$\ of the
original Cantor set $C$ is now endowed with a set of valued neighbours

\begin{equation}\label{eq8}
X_{\pm }^{i}=x \cdot x^{\mp \mid \tilde{x_i}\mid _{u}}.
\end{equation}

Finally, the element $x$\ is assigned the ultrametric valuation

\begin{equation}\label{eq9}
\parallel x\parallel \text{ }=\text{ }\underset{i}{\inf }\log _{x^{-1}}\frac{%
X_{+}^{i}}{x}\text{ }=\text{ }\underset{i}{\inf }\log _{x^{-1}}\frac{x}{%
X_{-}^{i}}
\end{equation}

\noindent so that $\parallel x\parallel $ $=2^{-n}$ $=3^{-ns}$ where $s=\frac{%
\log 2}{\log 3}$ , the Hausdorff dimension of the triadic Cantor set $C$\
and $n\rightarrow \infty $. As it turns out, this valuation exactly reproduces the 
nontrivial measure of \cite{pb} derived in the context of noncommutative geometry (c.f.,
definition of valued measure in Sec.3.2.) 

\vspace{.25cm}

Now, to make contact with the absolute value (\ref{eq5}) and the inversion rule (\ref{eq4}),
let for a sufficiently large but fixed $n$, $\tilde{x_i}$\ $\in \tilde I_{i}$ has
the form (we set for definiteness $m=0$)

\begin{equation}\label{eq10}
\tilde{x}_{i}=3^{-n} \cdot 3^{-n \cdot i2^{-r}}\text{ }\times a_i
\end{equation}

\noindent where $ni=2^{r} \cdot k_{i},$ $k_{i}$\ being, in general, a sufficiently large
real number\ and $a_{i}=\Sigma a_{ij}3^{-j}\in O_{i}$,\ a gap of size $3^{-r}
$of the Cantor set $C$\ and $a_{ij}\in \{0,1,2\}$. Then $0<\tilde{x_i}$\ \ $%
<3^{-n}$\ and $\mid \tilde{x_i}\mid $ $=$ $i.2^{-r}$. One also verifies that for an $a_i=1+a_{0i}$, $|a_i|_u=\underset{n\rightarrow \infty}{\lim}\log_{3^{n}}(a_i/3^{-n})=1$. By the inversion rule (\ref{eq4}) the elements \ $\tilde{x}_{i}$\ of the ball $\tilde I_{i}$\ now connect to an $x_{i}\in C$ \ given by

\begin{equation}\label{eq11}
x_{i}=3^{-n} \cdot 3^{-n(- \cdot i2^{-r})}\times b_i, \ b_i=\Sigma \text{ }b_{ij}\text{ }3^{-j},%
\text{ }b_{ij}\in \{0,2\}
\end{equation}

\noindent where $\lambda =a_i\times b_i \in (0,1).$ (Infinitesimal) Scales $\epsilon _{n}=3^{-n},$\ \
are the \emph{primary} scales when the scales $3^{-k_{i}}$(or equivalently 
$2^{-r}$) are the {\em secondary} scales.

\vspace{.25cm}

Finally, to verify the new multiplicative representation
one notes that there exists $y_{i}\in C$ in a neighbourhood of $x_{i}$ so that \bigskip

\begin{center}
$y_{i}=3^{-n}$ $c_{i},$ $\ \ c_{i}=$ $\Sigma $\ $c_{ij}$ $3^{-j},$ \ $%
c_{ij}\in \{0,2\}$\ 
\end{center}

\noindent and

\begin{equation}\label{eq12}
x_{i}=y_{i}.y_{i}^{-i.2^{-r}}
\end{equation}

\noindent so as to satisfy the identity

\begin{equation}\label{eq13}
c_{i}^{n-k_i}=b_{i}^{n}.
\end{equation}

To verify that (\ref{eq13}) is not empty we note that for the end points $\frac{1}{3}
$ and $\frac{2}{3}$, both belonging to $C$, (\ref{eq13}) means $\left( \frac{2}{3}%
\right) ^{n}=$\ \ $\left( \frac{1}{3}\right) ^{n-k_1}$\ yielding $k_1=ns,$\ $s=%
\frac{\log 2}{\log 3}$. For this value of $k_1$, (\ref{eq13}) now tells that $%
c_{i}^{1-s}=b_{i}$\ so that $c_{i}=({\frac{1}{3}})^{r}$\ and $b_{i}=\left( 
\frac{2}{3}\right) ^{r}$\ for a suitable $r$. Similar estimates for $k_i$
are available for other (consecutive)\ end points of (higher order) gaps.
It thus follows that the representation (\ref{eq8}) is realised  at the level of the finite Hausdorff measure of the set, when the value of the constant $k$ is real (rather than a natural number).
\vspace{.25cm}

Accordingly, it follows that a gap $O$ in $I/C$ ( which is a connected interval 
in the usual topology) containing a point $x$ of the Cantor set $C$ is indeed 
realised as an ``infinitesimal"  Cantor set in the valuation defined by the Cantor function associated with the original Cantor set $C$ itself. One thus concludes that 

\begin{proposition}
Any element $x$ of an ultrametric Cantor set $C$ has the multiplicative representation (\ref{eq8}) and the non-archimedean absolute value $||x||=\underset{i}{\inf} \  v(\tilde x_i)$. 
\end{proposition}
 
\subsection{Valued Measure}

Next it is easy to verify that the valued (metric) measure defined by the
ultrametric induced by the valued norm (\ref{eq9}) naturally yields the finite
Hausdorff measure of $C$.

\vspace{.25cm}

To recall, the valued measure $\mu _{v}:C\rightarrow R_{+}$\ is defined by

$(a)$\ $\mu _{v}$ $(\Phi )=0,$ \ $\Phi $ the null set

$(b)$ $\mu _{v}$ $[(0,x)]=$ $\parallel x\parallel ,$ when $x\in C$, 

$(c)$ For any $E\subset C,$ we have

\begin{equation}\label{eq14}
\mu _{v}(E)=\text{ }\underset{\delta \rightarrow 0}{\lim }\inf \text{ }%
\Sigma \left\{ \text{ }d_{u}(I_{i})\right\}
\end{equation}
\noindent where $I_{i}\in \tilde I_{\delta}$ and the infimum is over
all countable $\delta -$ covers of $E$ by clopen balls $I_i$ and 
$d_{u}(I_{i})$\ $=$ the non-archimedean ``diameter" of $I_{i}$ $=\sup \left\{ 
\text{ }\parallel x-y\parallel :x,\text{ }y\in I_{i}\right\}.$ 

\vspace{.25cm}

It follows that $\mu _{v}$\ is a metric (Lebesgue) outer measure on $C$
realised as an ultrametric space. Now for an infinitesimal (primary) scale $%
\epsilon _{n}=\beta^{n}$, we can choose a sufficiently small secondary scale $%
\tilde{\epsilon}_{m}=2^{-m}$, so that for any two $x$ and $y\in C$ with $%
\mid x-y\mid =d$ , we have $d_{u}(x,y)$\ $=$ $\parallel x-y\parallel $ $=$ $%
\tilde\epsilon _{m}^{s}\leq d^{s}\leq \left\{ d(x,y)\right\} ^{s}$, where $s=\frac{\log 2}{\log 1/\beta}$.  Accordingly,
it follows that $d_{u}(I_{i})\leq \left\{ \text{ }d(I_{i})\right\} ^{s}.$
Using this fact and also using the monotonicity of measures one can now
deduce that $\mu _{v}(E)=H^{s}(E)$. Choosing the full set $C$\ for $E$ and
noting that $s$ is the Hausdorff dimension it thus follows that $\mu
_{v}(C)=1.$\ We remark that the metric properties of the present 
ultrametric are indeed distinct from the
natural ultrametric  (c.f., \cite{sd1, sd2}), since the Lebesgue measure of $C$ in the natural ultrametric
is zero, but in the present case, the corresponding valued measure equals
the Hausdorff measure. More importantly, topologies induced by the two ultrametrics 
are also different, as seen in the following example.

\vspace{.25cm}

{\bf Example 3:} The sequence $\epsilon_n=\epsilon^{n-nl}, \ 0<\epsilon<1, \ 0<l<1$ converges to 0 in the usual metric (ultrametric), but converges to $l$ in the present ultrametric. For a sufficiently large $n$, choose $\epsilon^n$ as scale factor  and then relative infinitesimals  are $\tilde\epsilon_n= \lambda^{-1}\epsilon^{n+nl}, \ 0\ll\lambda<1$. Then, letting the secondary scale $\epsilon \rightarrow 0$, we have $|\tilde \epsilon_n|_u=\lim \log_{\epsilon^{-n}}(\epsilon^n/\tilde \epsilon_n) = l$ and hence $||\epsilon_n||=l$, by eq(\ref{eq9}), for a sufficiently large $n$. Thus, $\{\epsilon_n\} \rightarrow l$ in the ultrametric $||.||$. Letting $\epsilon={\tilde \epsilon}^{m}$, the sequence $\epsilon^{n-nl}$ is replaced by ${\tilde\epsilon}^{N-Nl}, \ N=nm $, so that the limit $\epsilon \rightarrow 0$ of the secondary scale is well defined, since it is realised as $m\rightarrow \infty$.
Note, however, that the sequence $\{\epsilon^n\}$ converges to 0, even in  $||.||$. For a sufficiently large but fixed $n$, we choose $\epsilon^{n+1}$ as the scale factor, so that $\tilde\epsilon_n=\lambda^{-1}\epsilon^{n+2},$ are relative infinitesimals and $|\tilde\epsilon_n|_u={\frac{1}{n+1}}$. More generally, for scales $\epsilon^{n+r}, \ r$ being a nonnegative real, we have  $|\tilde\epsilon_n|_u= {\frac{r}{n+r}}$. Thus, $||\epsilon^n||=\underset{r}{\inf}{\frac{r}{n+r}} = 0$. Further details of the convergence of  sequences and series (in a scale invariant real analysis) will be considered elsewhere.

\vspace{.25cm}

This example also gives an alternative proof that the metric $||.||$ is really an ultrametric.

\section{Differential Equation}

We now study the relationship of the scale free DE of the
form

\begin{equation}\label{eq15}
x\frac{dX}{dx}=X
\end{equation}

\noindent with a Cantor set $C$. As a preparation, let us recall how the simplest
Cauchy problem

\begin{equation}\label{eq16}
\frac{dX}{dx}=1,\text{ }X\text{ }(1)=1
\end{equation}

\noindent is solved on the interval $I$\ $=[0,1]$. One considers a partition $0=x_{0}<$
$x_{1}<x_{2}<\cdots <x_{n-1}<x_{n}=1.$ The desired result $X$ $(x)=x\underset{%
}{}$\ is obtained as a limit of a sum:$\underset{%
\begin{array}{c}
\bigtriangleup x_{j}\rightarrow 0 \\ 
n\rightarrow \infty%
\end{array}%
}{\text{ }\lim }\underset{j=1}{\overset{i}{\Sigma }}\bigtriangleup x_{j}$\
where $x_{i-1}<x<x_{i}.$ The scale free Cauchy problem 

\begin{equation}\label{sf}
x\frac{dX}{dx}=X,\text{ }X\text{ }(1)=1
\end{equation}

\noindent is also solved exactly in an analogous fashion.

\vspace{.25cm}

Let us now solve (\ref{sf}) in a slightly unconventional ``multiplicative"
approach \cite{dp2, dp3}. First we note that the neighbourhood of a point $x_{0}$\ is mapped
to that of $x=1$\ by a rescaling $x\rightarrow \frac{x}{x_{0}}.$ So we
concentrate only in the neighbour of $x=1.$ Let $x_{\pm }=1\pm \eta ,$ $%
X_{\pm }=X$\ $(x_{\pm }).$ Then the DE  in (\ref{sf}) splits into two
branches

\begin{equation}\label{eq18}
x_{\pm }\frac{dX_{\pm }}{dx_{\pm }}=X_{\pm }.
\end{equation}

Let us solve the left hand branch $X_{-}$. The standard solution is $%
X_{-}=1-\eta $. We now write

\begin{equation}\label{eq19}
X_{-}=\frac{1}{1+\eta }\text{ }X_{1-}
\end{equation}%

\noindent so that the correction factor $X_{1-}$ solves the renormalised
(self-similar) equation%
\begin{equation}\label{eq20}
x_{1-}\frac{dX_{1-}}{dx_{1-}}=X_{1-}
\end{equation}

\noindent where $x_{1-}$\ $=1-\eta ^{2}.$ Iterating ad infinitum, the desired solution
is retrieved (in a non-trivial form) thus

\begin{equation}\label{eq21}
X_{-}=\overset{\infty }{\underset{i=0}{\Pi }}\text{ }\frac{1}{1+\eta ^{2^{i}}}.
\end{equation}

The right hand branch, however, has the form 

\begin{equation}\label{eq22}
X_{+}=\frac{1}{1-\eta }\overset{\infty }{\underset{i=1}{\Pi }}\text{ }\frac{1%
}{1+\eta ^{2^{i}}}.
\end{equation}%
\ 

\bigskip The infinite product representation of $X_{-}$, for instance, is
interpreted as follows. The first iterated value exceeds the exact value by
an amount $\frac{\eta ^{2}}{1 + \eta }$ which is canceled progressively in
a self -similar manner over smaller and smaller inverted scales $\log $ \ $%
\left( 1-\eta ^{2^{i}}\right) ^{-1},$ $i=1,2,3,\cdots $.

\vspace{.25cm}

We note that the higher order correction factors $X_{c}=$\ \ $\overset{%
\infty }{\underset{i=1}{\Pi }}$ $\frac{1}{1+\eta ^{2^{i}}}$\ \ may therefore be
re-interpreted as a deletion process: viz; a portion of a line segment is
deleted progressively and self -similarly, analogous to the formation of a
Cantor set. Alternately, a product of the form $(1-\eta )$ $(1+\eta )=1-\eta
^{2}$\  could be considered to represent a deletion: a length of
size $\eta $\ in the neighbourhood of $1^{-},$\ is deleted progressively as $%
(1-\eta )(1+\eta )(1+\eta ^{2})\cdots \left( 1-\eta ^{2^{n-1}}\right)
=\left( 1-\eta ^{2^{i}}\right) ,$ $n\rightarrow \infty$.

\vspace{.25cm}

A second possibility is to interpret the multiplicative iteration process defined
above as a dynamical process in which the dynamic (independent) variable
undergoes increments not by the usual linear translations but by inversions
(hoppings) over smaller and smaller sizes. This then provides one with a
mechanism of deletion process stated above. Notice that the renormalised
corrections $X_{i}$ satisfy the self-similar scale free equations of the
from (\ref{eq20}) over smaller (logarithmic) scales (hopping sizes) $\log \left( 1-\eta
^{2^{i}}\right) ^{-1}$. Actually, these logarithmic scales inhabits
concomitant smaller scales of the from $\eta ^{i}$, that ordinarily arise in the 
context of an IFS. To justify this in a
greater detail, let us assume that the scale free problem (\ref{sf}) is now
defined on a closed subset $\tilde{\bf C}\subseteq I,$ called an {\em inverted Cantor
set}, where $\tilde{\bf C}=\underset{i}{\cup }$ $\tilde{I}_{i}$ is a countable
union of disjoint closed intervals $\tilde{I_i}$  of varying sizes. $%
\tilde{I_i}$  in fact, is the closure of a corresponding gap $O_{i}$ $($%
inclusive of the end points $)$ of the original Cantor set $C$. Suppose, a
dynamic variable (say, a particle) in motion on this set $\tilde{\bf C}$, hops
between the end points of each of such disjoint closed intervals $\tilde{I}%
_{i}$, following the scale free DE (\ref{sf}).
Let $\mid \tilde{I_0}\mid = \eta$   be the maximum hopping size. Then
the smaller hopping sizes are $\eta $ proportion of the remaining sizes of
the set $\tilde{\bf C}$ viz.,
$\mid \tilde{I_i}$  $\mid =\eta $ $(1-\eta )^{i-1}\equiv \eta_i$ (say), when we assume $%
\mid \tilde{\bf C}\mid $ $=1.$\ Because of the rescaling symmetry (scale
invariance) of the DE (\ref{sf}) each of the component
intervals $\tilde{I}_{i}$ could be imagined to have been symmetrically
placed at 1 with end points at $x_{i\pm }=\frac{1}{2}$ $(1\pm \eta_i).$ Now,
the particle at left end point $x_{0-\text{ }}$ of $\tilde{I_0}$ (say)  hops
to the right end point $x_{0+\text{ }}$ following the rule (c.f., (\ref{eq19}))

\begin{equation}\label{eq23}
x_{0-\text{ }}\rightarrow x_{0-\text{ }}^{-1}=x_{0+\text{ }}{X}_{1}
\end{equation}

\noindent so that we have, using (\ref{eq22})

\begin{equation}\label{eq24}
(1-\eta )^{-1}=(1+\eta )\text{ }{X}_{c},\text{ } X_{c}=\underset{i=1}{%
\overset{\infty }{\Pi }}\left( 1+\eta ^{2^{i}}\right)
\end{equation}%
\noindent because of the scale invariance. Eqn(\ref{eq24}) tells that hopping motion of the type considered above, of any given size $\eta $ is accomplished by an infinite
cascade of self similar smaller\ scale inverted motions of  sizes $\eta
^{2^{i}},$ $i=1,2,\cdots $. The total length covered by all these
self-similar jumps , viz., 1, is reached multiplicatively i.e. as $1=\underset{%
n\rightarrow \infty }{\lim }\left( 1-\eta ^{2^{i}}\right) $, reminiscent of
an ultrametric limiting process. Notice that in the ordinary sense, the
total jump size is determined additively as an infinite series viz., $%
\underset{1}{\overset{\infty }{\Sigma }}$ $\eta $ $(1-\eta )^{i-1}= 
\Sigma \eta_i \ =\ 1.$

\vspace{.25cm}

The above multiplicative model of the Cauchy problem (\ref{sf}) should have a
natural relevance in the context of a Cantor set (or, equivalently to an
ultrametric  space). Let us now assume that the problem (\ref{sf}) is
defined on a Cantor set $C$, .i.e. $x\in C$. This Cantor set is now realised as an
inequivalent ultrametric space. Accordingly, each $x\in C$\  is replaced by an
infinitesimal copy of the inverted Cantor set $\tilde{\bf C}$. {\em Because of scale
invariance, the DE (\ref{sf}) at an $x\in C$, which is actually not defined in the usual 
(even in the natural ultrametric) sense,  
is now raised to an equation which is well defined on a closed set of the form $\tilde{\bf C}$ and treated as a multiplicative model.}

\vspace{.25cm}

We note that the solution (\ref{eq21}) and (\ref{eq22}) is the standard solution derived in an unconventional way and interpreted non-trivially. On a Cantor set,
however, the equation (\ref{sf}) can accommodate a host of new solutions in
consonance with the multiplicative model interpretation. The origin of
these new solutions could be explained in the context of locally constant
functions (LCF). To justify, in a most natural way, the existence of  locally
constant functions, let us write a solution of (\ref{sf}) in the form

\begin{equation}\label{eq25}
X=x \cdot x^{\tau (x^{-1})}.
\end{equation}

\noindent The function $\tau (x)$\ here represents a LCF and is
defined by the scale free equation on logarithmic variables, viz:

\bigskip 
\begin{equation}\label{eq26}
\log x^{-1}\frac{d\tau }{d\log x^{-1}}=\tau.
\end{equation}

\noindent Clearly $\tau (x)$\ corresponds to our non-trivial valuation (\ref{eq5}) denoted $v(\tilde{x})=\mid \tilde{x}$\ $\mid _{u}$. To verify $v(\tilde x)$, indeed is a LCF, we note that

\begin{equation}\label{eq27}
\frac{d}{dx}v(\tilde{x})=\underset{\epsilon \rightarrow 0}{\lim }\frac{d}{dx}%
\left( \frac{\log x}{\log \epsilon }+1\right) =0.
\end{equation}

\noindent Equation (\ref{eq26}), on the other hand, reveals the variability of a LCF over smaller logarithmic scales. Of course, the valuation
also passes this test \bigskip
\begin{center}
$\log $ $v(\tilde{x})=\log \log \frac{x}{\epsilon }+\log \log _{\epsilon
}\lambda -\log \log \epsilon ^{-1}$
\end{center}
\bigskip 
\noindent leading to equation (\ref{eq26}) in the inverted rescaled real variable\ $%
\frac{x}{\epsilon }$\ (in the $\log \log $\ scale). We have already seen
that $v(\tilde{x})$\ relates to an appropriate Cantor function.
Consequently, a Cantor function $\phi (x)$ is shown to be a LCF with variability over $\log \log $\ scales. Equation (\ref{eq25})
constitutes an ultrametric extension not only of a Cantor set, but of any
connected interval of $R$ \cite{ad}. This is already verified explicitly in the $\frac{
1}{3}$\ rd Cantor set, in Sec.3.1.   

\vspace{.25cm}

The main results derived in this section are summarised in 

\begin{theorem}
An element of an ultrametric  Cantor set $C$ is replaced by the set of gaps of the Cantor set $\tilde C$ where relative infinitesimals are supposed to live in. Increments on such an extended Cantor set $C$ is accomplished by following an inversion rule of the form (\ref{eq4}). A scale free differential equation of the form eq(\ref{sf}) is well defined on such an ultrametric space and accommodates Cantor functions as locally constant functions.
The associated infinitesimal valuation $v(\tilde x)$ is a locally constant function with variability over double logarithmic scales.
\end{theorem}

\section{{Growth of Measure}}

In Sec 3 we studied the valued ultrametric structure of a measure zero
Cantor set. Here we study a few more general properties of the valued
ultrametricity. We note, at first, that the valued structures of $\tilde{x}$
and\ $x$ of the middle third Cantor set $C_{1/3}$, viz., equation (\ref{eq10})
and (\ref{eq11})  actually correspond to the solution (\ref{eq24}) when $a_i$, $b_i$ and $\lambda$ are identified with  $(1-\eta ),$ $(1+\eta )$
and $\tilde{X}$\ respectively. The function $\tilde{X}$\ \ is a LCF. The valuation $v(\tilde{x})$ is identified with the LC Cantor function corresponding to $C_{1/3}$. However, as verified in equation (\ref{eq27}), $v(\tilde{x})$ is indeed a LCF satisfying

\begin{equation*}
\frac{d}{dx}v(\tilde{x}(x))=0.
\end{equation*}%
Consequently, $v$\ is, not only a LCF,  but more
importantly is a {\em reparametrisation} invariant object. As a result, $v$\ does 
not require to be an explicit function of the original variable $x$ but may be a function
instead of {\em any} monotonic, continuously first differentiable   function of $x$. By the same token, $v$
does not depend explicitly on the scale $\epsilon $\ inherited from the original
(mother) Cantor set, as we did in the example of $C_{1/3}$. In the
following example, we show that relative infinitesimals may instead live in a
positive measure Cantor set. Notice that in the general representation of
the valued ultrametric in (5), the parameter may be a constant
independent of an explicit $\epsilon.$

{\bf Example 4:}
Suppose that eq(\ref{eq10}) and (\ref{eq11}) are replaced by

\begin{equation}\label{eq28}
\tilde{x}=\beta ^{n}\left( \beta ^{n^{\delta }}\right) ^{n\beta _{n}\left(
1+\gamma _{m}\right) }\times a
\end{equation}

\noindent where $a=(1-\beta )\left( 1+\underset{1}{\overset{\infty }{\Sigma }}%
a_{i}\beta ^{i}\right) ,$ $a_{i}\in \{0,1\},$ $\beta =\frac{1}{2}(1-\alpha )$,%
\ and $\beta _{n\text{ }}$and $\gamma _{n}$ are two non-increasing sequence
of positive numbers such that $\beta _{n\text{ }}\rightarrow 0$\ as $%
n\rightarrow \infty $ and $m$ may be independent of $n$ or may vary with $n$
more slowly, and  $\delta >1$ is a constant.

\vspace{.25cm}

Although $\beta _{n\text{ }}\rightarrow 0$\ as $n\rightarrow \infty $ , the
valuation $v(\tilde{x})$ could be non-trivial, since

\begin{eqnarray}\label{eq29}
v(\tilde{x}) &=&\underset{n\rightarrow \infty }{\lim }\log_{\beta
^{-n}}{\frac{\beta ^{n}}{\tilde{x}}} \\
&=&\underset{n\rightarrow \infty }{\lim }\left[ n^{\delta }\beta _{n}\left(
1+\gamma _{m}\right) +\log \underset{\beta ^{-n}}{a}\right]  \notag \\
&=&l+\tilde{\gamma}_{m_{n}}(\delta )  \notag
\end{eqnarray}

\noindent when we assume $n^{\delta }\beta _{n}\rightarrow l$ as $%
n\rightarrow \infty $ and $n^{\delta }\beta _{n}\gamma _{m}\rightarrow $\ $%
\tilde{\gamma}_{m_{n}}(\delta )$\ is a sub-dominant slowly varying
non-increasing sequence, for a real $m_n>0$. The representation (\ref{eq28}) tells that a scale free infinitesimal $\frac{\tilde{x}}{\beta ^{n}}$ may live in a Cantor set $\tilde C_{p}$ of Example 1 (Sec 2), so that $m(\tilde C_{p})=l$. Let the original Cantor set be a
middle $\alpha $ set $C_{\alpha }$ with the uniform scale factor $\beta =%
\frac{1}{2}(1-\alpha )$. For the positive measure set $\tilde C_{p}$ the scale
factor at the nth iteration is $\tilde{\beta}_{n}=2^{-n}\underset{i=1}{%
\overset{n }{\Sigma }}\left( 1-\alpha _{i}\right) $\ \ \ \ \ \ \ \ \
and $l=m\left( \tilde C_{p}\right) =$\ \ \ $\underset{i=1}{\overset{\infty }{\prod
}}\left( 1-\alpha _{i}\right) =\underset{n\rightarrow \infty }{\lim }$\ \ $%
2^n\tilde{\beta}_{n}.$  Let us choose\ $\delta >1$\ such that $\tilde{\beta}_{n}=$%
\ $\beta ^{n^{\delta }}$. Then $n^{\delta }\beta _{n}\rightarrow l$ tells
that $\beta _{n}\approx \frac{l\log \beta  }{\log \tilde \beta_n}$\ as $%
n\rightarrow \infty $. Thus the dominant term $l$ of the valuation $v(%
\tilde{x})$ is a constant while the subdominant asymptotic $\tilde{\gamma}%
_{m}(\delta )$\ \ could be a genuine LCF (i.e. a
Cantor function for a sub dominant Cantor like set $C_s$ (say) ), precise
determination of which depends on the explicit model of the Cantor set $\tilde C_{p}$%
. It follows, therefore, from (\ref{eq8}) and (\ref{eq9}) the ultrametric valuation of $x\in
C_{\alpha }$\ now has the form

\begin{equation}\label{eq30}
\parallel x\parallel =l+\tilde{\gamma}_{m_{n}}(\delta ). 
\end{equation}

\noindent For larger and larger values of $n$ $(\rightarrow \infty )$, we can disregard
the sub-dominant term (since $\tilde{\gamma}_{m_{n}}$\ $\rightarrow 0$ as $%
m_n \rightarrow \infty $) so that 

\begin{equation}\label{eq31}
\parallel x\parallel =l \ \text{ }\forall \text{ }x\in C_{\alpha },\text{ }%
x\neq 0.
\end{equation}

\noindent Clearly the trivial ultrametric (\ref{eq31}) reveals that the mother set $C_{\alpha} $ must get  deformed to a positive measure set $C_p$ so that $\mu _{v}(C_p)=m(C_p)=l$, 
when the reparametrisation invariance of LC correction factors is invoked. 
Indeed, we have \ $\parallel x-y\parallel =l$ for any two $x,y\in C_p.$ Thus,
any single clopen ball $B\left( x_{0}\right) ,$ $x_{0}\in C_p$ (say) covers
the compact $C_p$ and hence $\mu _{v}(C_p)=$ $d_{u}\left( B\left( x_{0}\right)
\right) =l.$  

\vspace{.25cm}

To summarise, we have shown that any element $x\in C_{\alpha }$\ when
deformed by the non-trivial, reparametrisation invariant valuation of
relative infinitesimals, is identified with an element of a 1-set $C_p$. Because of this invariance, the relative infinitesimals may be assumed to live in a positive measure set $\tilde C_p$, which, in turn, determines the measure (size) of the deformed set $C_p$. Since
each element $x\in C\subset [0,1]$ is written as the arithmetic sum of two elements $
x_{0}\in C_{\alpha }$ and $x_{1}\in $\ $C_{\alpha ^{\prime }}$\ ( $C_{\alpha
^{\prime }}$ being the Cantor set of infinitesimal neighbours of $x_{0}$), it
follows from a theorem of Solomyak \cite{solo} that for $\beta =\frac{1}{2}\left(
1-\alpha \right) \in \left( 0,\frac{1}{2}\right) ,$ there exists $C_{\alpha
^{\prime }}$ for a.e. $\beta ^{\prime }=\frac{1}{2}\left( 1-\alpha ^{\prime
}\right) \in \left( 0,\frac{1}{2}\right) $ so that $C_{\alpha }+C_{\alpha
^{\prime }},$ has positive measure and $\frac{1}{\log \frac{1}{\beta }}+%
\frac{1}{\log \frac{1}{\beta ^{\prime }}}>\frac{1}{\log 2}.$ This, therefore, constitutes an alternative proof for the said assertion. Indeed, in the above construction, the set of infinitesimals $C_{\alpha^{\prime}}$ itself is a 1-set $\tilde C_{p}$.

\vspace{.25cm}

It follows, accordingly, that a slower rate of removal of middle open sets 
compared to a measure zero Cantor set hides a positive measure in an 
infinitesimal scaling factor which is exposed under the present scale invariant 
valuation. The uniform rate of deletion in the case of a measure zero set 
is  violated because of the underlying reparametrisation invariance. Further, in a 
dynamical process leading to  a Cantor set, a positive measure Cantor set $C_p$
is favoured a.s (almost surely) compared to a measure zero set $C_{\alpha}$ since 
relative infinitesimal neighbours a.s. lie in a Cantor set $C_{\alpha^{\prime}}$ 
satisfying the above constraints.

\vspace{.25cm}

The generic result that follows from this example is stated thus

\begin{theorem}
Because of the reparametrisation invariance of the infinitesimal valuation, a measure zero Cantor set $C_{\alpha}$ is a.s. deformed to a positive measure Cantor set $C_p$, the measure of which is determined by the Cantor set $\tilde C_p$ in which the relative infinitesimals are supposed to live. 
\end{theorem}

Next to expose the significance of the sub-dominant term, let us first define a
renormalized valuation $v_R(\tilde x)$:

\begin{equation}\label{eq32}
v_{R}\left( \tilde{x}\right) =\log _{\beta ^{n}}\log _{_{\beta ^{n}}}\left[ 
\frac{\tilde{x}}{\left( {\beta ^{n}}\right) ^{1+v_{0}\left( \tilde{x}%
\right) }}\right], \ \ n\rightarrow \infty
\end{equation}

\noindent where $v_{0}\left( \tilde{x}\right) =l<1$ is the dominant valuation of the
infinitesimal $\tilde{x}.$ The LCF $\tilde{\gamma}_{m_{n}}(\delta )\ $ is now given by (c.f.(6))

\begin{equation}\label{eq33}
\tilde{\gamma}_{m_{n}}(\delta )\ =\alpha _{i}\beta ^{m_{n}\rho \left( \delta
\right) }
\end{equation}

\noindent where the $\delta -$dependent constant $\rho $ is called a {\em renormalised valuated exponent} and the non-zero constant $\alpha _{i}$ assumes values from a finite
set for a secondary scale $\beta^{m_n}$. As will become clear the valuated exponent $\rho $ is useful to distinguish two sets with identical Hausdorff dimensions.

\vspace{.25cm}

\subsection{Applications}

{\bf 1. Middle third Cantor sets:} As an application of the renormalised valuated exponent, let us first consider 
a class of $s$- sets where $s=\log_3 2$, constructed as a slight variation of the process of Example 2, Sec.2. Let $I=[0,1]$. Also let $0<<\delta_n=3^{-(n+1)\alpha_n}
\lessapprox 1,  \ n,=1,2,\cdots$ (so that $\delta_n^{-1}\gtrapprox 1$), be a non-increasing sequence. For definiteness, one may choose $\alpha_n= q^{-n}$, for a sufficiently large positive integer $n$ and $q>1$. In that case $\alpha_n$ may be considered to belong to the range set of an appropriate Cantor function. Delete the middle open interval of length 1/3. Next, delete a length $3^{(-2(1+\alpha_1))}$ from each of the two closed subintervals. Then, delete the length  $3^{(-3(1-\alpha_2))}$ from each of $2^2$ closed subintervals. Call these two operations together $O_1$. $O_n$ consists of two steps: deletion of $2^{n+1}$ open intervals of length $3^{-(n+1)(1+\alpha_{n})}$, which is succeeded by the next deletion of lengths  $3^{-(n+2)(1-\alpha_{n+1})}$ from  $2^{n+2}$ remaining closed subintervals. Notice that we are considering a set of fluctuating scale factors, i.e., in the $(n+1)$th step open intervals of slightly smaller sizes compared to the middle third set are removed. In the next step, however, open intervals of slightly bigger sizes are removed. As a consequence,  we get a family of limit sets which are indistinguishable and equivalent to the middle third Cantor set at the level of the Hausdorff  dimension, but nevertheless, distinguishable at the level of renormalised valuated exponents. Indeed, the total length of deleted open intervals viz.,  ${\frac{1}{3}}+{\frac{2}{3^2}}\delta_1+ {\frac{2^2}{3^3}}\delta^{-1}_2 +\cdots= 1+\sum u_n$ equals 1, when the series of real numbers $\sum u_n$ vanishes. The sequence $u_n$ is determined by the sequence $\alpha_n$, i.e., $\alpha_n=\log_{3^{(n+1)}} (1-{\frac{3^{(n+1)}}{2^n}} \tilde u_n)^{-1}$ and $\alpha_{n+1} =\log_{3^{(n+2)}} (1+{\frac{3^{(n+2)}}{2^{n+1}}} \tilde u_{n+1})$ so that $u_n=-\tilde u_n, \ u_{n+1}=\tilde u_{n+1}$. Clearly, such a series exists. Hence, all such sets are of measure zero.

\vspace{.25cm}

Now, to determine the Hausdorff dimension, we first note that the scaling of closed intervals (bridges) follows the recurrence $2l_n=l_{n-1}-\delta_n^{\pm 1}3^{-(n+1)}$, where $+$ sign goes with an odd $n$ and the - sign with  $n$ even and $l_n$ denotes the length of each closed interval at level $n$. Accordingly, $l_n={\frac{1}{3.2^n}} [1- {\frac{\delta_1}{3}} - {\frac{2\delta^{-1}}{3^2}} - \cdots {\frac{2^{n-1}\delta_n^{\pm 1}}{3^n}}] \approx {\frac{\delta_{n+1}^{\pm 1}}{3^{n+2}}}$, for a sufficiently large $n$. As a consequence, the scale factors behave as either $\beta_{n+1}=3^{-(n+1)(1+\alpha_{n})}$  or $\beta_{n+2}=3^{-(n+2)(1-\alpha_{n+1})}$ respectively, and hence, the lower and upper box dimensions and the Hausdorff dimension are all equal and equal to  $\underset{n\rightarrow \infty }{\lim } {\frac{\log 2}{ \log 3^{(1 \pm\alpha_{n})}}}=\log_3 2$.

\vspace{.25cm}

One may also estimate the thickness of these sets easily. Because of the above scaling, the limiting length of the closed intervals (bridges) coincides with that of the corresponding  gap (viz., $\delta_{n+1}^{\pm 1}3^{-(n+2)}$) at the $n$th level. It follows therefore that the ratio of sizes of bridges and gaps (c.f., Sec.2.1) has the limiting value 1. Hence, thickness of all these sets coincides with that of the classical middle third Cantor set as well. 

\vspace{.25cm}

However, a higher order (renormalised) valuated exponent can indeed reveal the local dissimilarities of such an $s-$set. Extending the representations (\ref{eq11}) and (\ref{eq28}) a little further to suit the present problem, we would now have for an element $x$ of the $s-$set,

\begin{equation}\label{eq34}
x_{i\pm}= 3^{-n}\cdot 3^{-n(-i2^{-m_n(1\pm \alpha_{m_n})})}\times b, \ ||b||=1
\end{equation} 

\noindent where $i$ assumes values from a finite set and $m_n\rightarrow \infty$ at a slower rate as $n\rightarrow \infty$,  so that a renormalised valuation is defined as

\begin{equation}\label{eq35}
v_R(x)=\underset{i}{\inf}\log_{2^{-m_n}}\log_{3^{-n}} (x_{i+}/x_0) =\alpha_{m_n}, \ x_0= 3^{-n} \cdot 3^{-n(-i2^{-m_n})}. 
\end{equation}

\noindent It now follows from the definition of $\alpha_{m_n}$, that one can  find a sufficiently large natural number $q>>1$ such that $\alpha_{m_n}=q^{-m_n}$. Consequently, we obtain $v_R(x)=\alpha_{m_n}= 3^{-\rho\tilde m_n}$, where $\rho=\log_{3^r} q \ >0$ is the {\em valuated exponent}, for suitable positive integers $r$ and $\tilde m_n$. 

\vspace{.25cm}

Now, to justify the existence of such a $q$, let us first assume $\tilde u_{2m}=u^1_m$ and $\tilde u_{2m+1}= u^2_m$ such that $\sum u^i_m=l$. Then $\underset{2}{\overset {\infty}{\sum}} u_n = \sum \tilde u_{2m+1}-\sum \tilde u_{2m}=l-l=0.$ Consequently, 
$\alpha_{2m}=\log_{3^{(2m+1)}} (1-{\frac{3^{(2m+1)}}{2^{2m}}} u^1_m)^{-1}$ and $\alpha_{2m+1} =\log_{3^{(2m+2)}} (1+{\frac{3^{(2m+2)}}{2^{2m+1}}} u^2_m)$. 

\vspace{.3cm}

Let $\eta_{2m}={\frac{3^{(2m+1)}}{2^{2m}}} u^1_m$ and $\tilde \eta_{2m+1} = {\frac{3^{(2m+2)}}{2^{2m+1}}} u^2_m$. Then the functions $(1-\eta_{2m})^{-1}$ and $1+\tilde \eta_{2m+1}$ are identified as LCF of the form eq(23) and eq(24) (Sec.4), in the neighbourhood of 1. Using scale invariance, we can choose for $x$ in eq(23) as $x=3^{-n}(1-\eta_{n})$ (or $x=3^{-n}(1+\tilde\eta_{n})$), and the scale factor $\epsilon=3^{-n}$. Thus, there exists a Cantor function $\tau(\tilde x), \ \tilde x=x/\epsilon$ such that $\log_{\epsilon^{-1}} \tilde x^{-1}=\tau(\tilde x) \ ( {\rm or} \ \log_{\epsilon^{-1}} \tilde x=\tau(\tilde x))$. As a result, there exists positive integers $q$ and $m_n$ so that the sequence $\{q^{-m_n}\}\subset {\rm Range} (\tau(\tilde x))$. More generally, because of the local constancy, the limiting form $\alpha_n$ could be $\alpha_n = \tilde l + q^{-m_n}$, where $\tilde l$ is a non-negative constant, $0\leq\tilde l <1$.

\vspace{.25cm}

We remark that the exponent $\rho$ may be considered to be the inverse of the Hausdorff dimension of a {\em residual} Cantor set that would remain attached with infinitesimal scales in a neighbourhood of a point (of the original Cantor set). For the classical middle third Cantor set $\alpha_n=0\ \forall n$ and so $\rho=\infty$, which is consistent with the fact that the residual set is null. Since, sets with infinite Hausdorff dimension $s=\infty$ are excluded, by definition, $\rho$ indeed is positive $\rho>0$.

\bigskip
{\bf 2. 1-sets:}  Irregular 1-sets \cite{falc}  are positive measure Cantor sets and are generally classified on the basis of fatness and/or uncertainty exponents. The LC renormalised valuation (\ref{eq32}) and (\ref{eq33}) now tells that $v_R(\tilde x)$ is a Cantor function corresponding to a subdominant residual Cantor set $C_s$, and so has the form $v_R(\tilde x)= \alpha _{i}\beta ^{m_{n}\rho}$.  As for the $s$-sets, the valuated exponent $\rho >0$ equals the inverse of the Hausdorff dimension of the residual set $C_s$. For $\rho=\infty$, the double exponential factor in (\ref{eq28}) drops out (i.e., reduces to the trivial factor $\beta^n$), and hence the 1-set is a regular set \cite{falc} having connected components (actually corresponds to a nonfractal set). Consequently, $0< \rho\leq\infty$.

\vspace{.25cm}

Now, to compare with the fatness exponent \cite{far1, far2}, we first recall the relationship between the uncertainty exponent $\alpha, \ 0<\alpha\leq 1$ \cite{ott} and the fatness exponent $\tilde\beta, \ 0<\tilde \beta\leq \infty$. It is shown \cite{far2} that $\tilde \beta=\alpha$ in [0,1], so that there is essentially the fatness exponent that has to be considered. We claim that $\rho=\tilde\beta$. The parameter $\tilde\beta$ is defined as 

\begin{equation}\label{eq36}
\tilde\beta=\underset{\epsilon\rightarrow 0}{\lim}{\frac{\log [\mu(\epsilon)-\mu(0)]}{\log\epsilon}}
\end{equation}

\noindent where $\mu(\epsilon)$ is a LC measure which tells the scaling of smaller gap sizes when the smaller gaps are coarse grained by fattening by the amount $\epsilon$ and $\mu(0)$ equals the positive (Lebesgue) measure of the set. In our multiplicative representation (c.f. (\ref{eq28}) and (\ref{eq33})), the fattening size is $\epsilon=\beta^n$ and 

\begin{equation}\label{eq37}
x=\beta^{n}\left( \beta ^{n}\right) ^{-(l+k\beta^{n\rho}) }\times b
\end{equation}

\noindent where $k$ is a constant independent of $\beta$, so that the exponent $\rho$ is defined by (\ref{eq36}) when we identify $\mu(\beta^n)=\log_{\beta^{n}} (x/\beta^n)$. Notice that $\mu(0)=\underset{n\rightarrow \infty}{\lim} \log_{\beta^{n}} (x/\beta^n)=l$. Notice also that the measure $\mu$ here is nothing but the valuation of relative infinitesimals at the fattened scale $\epsilon$, which equals the full measure of the Cantor set $\tilde C_p$ at the scale $\epsilon$ (c.f. Example 4) where the infinitesimals live. Because of the reparametrisation invariance, we may  suppose that $\tilde C_p$ is determined by the original 1-set and vice versa. At the scale $\epsilon$, the gaps of $\tilde C_p$ are fattened by  the amount $\epsilon$, and in the presence of a positive measure, the said valuation is determined by the sum of the fattened gap sizes. For a zero measure set, this valuation, on the other hand, is determined instead by the finite Hausdorff measure, upto a finer (double logarithmic) scale  correction that arises from the possible presence of local fine structures (c.f., above application). This observation proves the claim. 

\section{Final Remarks}

The general representation of the fine structure induced by relative infinitesimals to any point $x$ of a Cantor set $C$, with a nonnegative measure $m( C)\geq 0$, that has emerged from the above analysis, has the form 

\begin{equation}\label{gr}
X_{\pm}(x)=x\left( x^{\pm v(\tilde x)(1-H^s(\tilde C))}\cdot x^{\pm x^{\pm v_R(\tilde x)(1-H^{{\tilde s}}(\tilde C_r))}}\cdots\right )\times a.
\end{equation}

\noindent Here, relative infinitesimals $\tilde x$ are supposed to live in a Cantor set $\tilde C$. The actual nature of $\tilde C$ is determined by $C$ and vice versa. The first (logarithmic) order valuation $v(\tilde x)=m(\tilde C)$ when $m(\tilde C)>0$ (and $H^s(\tilde C)=1$ for $s=1$, thus eliminating the first order corrected factor). When $m(\tilde C)=0$, $v(\tilde x)>0$ and $H^s(\tilde C)$ is the finite $s$ dimensional Hausdorff measure of the (Lebesgue) measure zero set $\tilde C$. The higher order (double logarithmic) valuation has two interpretations: for a positive Lebesgue measure  set, $v_R$ relates to the fatness exponent, which is interpreted as the inverse of the Hausdorff dimension of a residual Cantor set $\tilde C_r$. On the other hand, when the Lebesgue measure is zero, $V_R$ gives rise to a valuated exponent $\rho$ which can be used to distinguish between sets having identical Hausdorff dimension $s$ and thickness. $\tilde s=1/\rho$ is the Hausdorff dimension of $ \tilde C_r$ and $H^{\tilde s}(\tilde C_r)$ is the corresponding  $\tilde s$ dimensional Hausdorff measure. The nontrivial factors involving Hausdorff measures may be seen to arise analogous to (\ref{eq29}), in view of the reparametrisation invariance of the higher order valuations $v_R$. Finally, the ellipses indicate the possible existence of higher order exponential factors. The ordinary (real analysis) representation is reproduced for $\tilde C=\Phi$.

A detail study of the {\em generalized Euler's factorisation } (\ref{gr}) will be taken up separately. Let us  only remark here that the nontrivial measure dependent exponents tell that the asymptotic limit $x\rightarrow 0^+$ actually depends on the {\em size} and {\em nature} of the underlying set on which the variable $x$  is supposed to live in. Exploiting the reparametrisation invariance, the structure of the underlying set $C$ could also be altered significantly, that is, a measure zero set may acquire a positive measure.

\section*{Acknowledgement}

It is a pleasure to thank the anonymous referees for very constructive recommendations to improve the quality of the paper.

\end{document}